\documentclass[12pt]{article}
\usepackage{amsmath,amssymb,amsthm}

\newcommand{\Var}{{\cal{V}_{\mathbb{C}}}}
\newcommand{\fgs}{{\mbox{\rm f.$G$-s.}}}

\def\1{\underline{1}}

\def\P{\mathbb P}

\def\Z{{\mathbb Z}}
\def\Q{{\mathbb Q}}
\def\C{{\mathbb C}}

\def\DD{{\cal D}}

\newtheorem*{theorem*}{Theorem}

\newenvironment{definition}
{\smallskip\noindent{\bf Definition\/}:}{\smallskip\par}

\newenvironment{examples}
{\smallskip\noindent{\bf Examples\/}.}{\smallskip\par}
\newenvironment{remark}
{\smallskip\noindent{\bf Remark\/}.}{\smallskip\par}

\title{An equivariant version of the monodromy \\ zeta function
}
\author{S.M.~Gusein-Zade \thanks{Partially supported by the grants
RFBR-007-00593, INTAS-05-7805, NWO-RFBR 047.011.2004.026, and
RFBR-JSPS 06-01-91063.
Address: Moscow State University, Faculty
of Mathematics and Mechanics, Moscow, 119992, Russia. E-mail:
sabir\symbol{'100}mccme.ru} \and I.~Luengo \and
A.~Melle--Hern\'andez \thanks{The last two authors were partially
supported by the grant MTM2007-67908-C02-02. Address: University
Complutense de Madrid, Dept. of Algebra, Madrid, 28040, Spain.
E-mail: iluengo\symbol{'100}mat.ucm.es,
amelle\symbol{'100}mat.ucm.es}}
\date{}
\begin{document}
\def\eps{\varepsilon}

\maketitle

\centerline{\it To S.P. Novikov with admiration.}

\begin{abstract}
We offer an equivariant version of the classical monodromy zeta function of a singularity
as a series with coefficients from the Gro\-then\-dieck ring of finite $G$-sets tensored by
the field of rational numbers.
Main two ingredients of the definition are equivariant Lefschetz numbers and the $\lambda$-structure on the Grothendieck ring of finite $G$-sets. We give an A'Campo type formula for the equivariant zeta function.
\end{abstract}

A number of invariants of singularities have equivariant versions for singularities invariant with respect to an action of a finite group $G$. For example, in \cite{CTC}, an {equivariant} version of the {Milnor number} of a $G$-invariant function germ is an element of the ring of virtual representations of the group $G$. However, one can say that it is not clear what should be considered as the equivariant version of the classical {monodromy zeta function} of a singularity.

Here we offer a candidate for this role. Main two ingredients of the definition are {equivariant Lefschetz numbers} (defined, e.g., in \cite{LR},\cite{GW}) and the $\lambda$-structure on the Grothendieck ring $K_0(\fgs)$ of finite $G$-sets (also called Burnside ring, see \cite{knutsen})
which is connected with the power structure over the Grothendieck ring of complex quasi-projective varieties discussed in \cite{GLM1}. Moreover, we consider the equivariant Milnor number of a singularity of a $G$-invariant function germ as an element of the Grothendieck ring of finite $G$-sets.
This notion reduces to the one considered in \cite{CTC} under the natural homomorphism from 
$K_0(\fgs)$  to the ring $R(G)$ of virtual representations of the group $G$, but in some cases is a certain refinement of that one.

The classical monodromy zeta function $\zeta_f(t)$ of the germ of a function $f$ on a variety is a rational function in a variable $t$ or a power series in $t$ (with integer coefficients).
The equivariant version $\zeta^G_f(t)$ offered here is a power series in one variable $t$ with coefficients from the tensor product $K_0(\fgs)\otimes \Q$ of the Grothendieck ring of finite $G$-sets with the field $\Q$ of rational numbers. 

We write down an A'Campo type formula for the equivariant zeta function $\zeta^G_f(t)$.

\section{The Grothendieck ring of finite $G$-sets}\label{uno}

A {finite $G$-set} is a finite set with an action of a finite group $G$.
An {isomorphism} between two $G$ sets is a one-to-one correspondence from one to the other which respects the $G$-action. 
Isomorphism classes 
of {irreducible} $G$-sets, also called simple or transitive,
(i.e. those which consist of exactly one orbit) are in one-to-one correspondence with the set $\mbox{\rm consub}(G)$ of conjugacy classes of subgroups of $G$, see \cite{knutsen}.

The Grothendieck ring of finite $G$-sets is the group generated by isomorphism classes of finite $G$-sets with the relation $[A\coprod B]=[A]+[B]$ for finite $G$-sets $A$ and $B$ and with  the multiplication defined by the cartesian product.
As an abelian group $K_0(\fgs)$  is freely generated by isomorphism classes of irreducible $G$-sets and, therefore, an element of $K_0(\fgs)$ can be written as
$\sum_{h\in {\rm consub}(G)} k_h [G/H]$, where $H$ is a representative of the conjugacy class
$h$, $k_h\in \Z$.

The Grothendieck ring $K_0(\fgs)$ has a natural $\lambda$-structure defined by the series 
$$
(1-t)^{-[X]}:=1+[X]\,t+[S^2X]\,t^2+[S^3X]\,t^3+\ldots
$$ 
where $S^kX=X^k/S_k$ is the $k$-th symmetric power of the $G$-set $X$ with the natural $G$-action. (This $\lambda$-structure induces a power structure over the Grothendieck ring $K_0(\fgs)$ connected with the power structure over the Grothendieck ring $K_0(\Var)$ of quasi-projective algebraic varieties introduced in \cite{GLM1}).

\medskip

\begin{examples} 
\begin{enumerate}
\item Let $G$ be the cyclic group $\Z_6$ of order $6$. There are four conjugacy classes of subgroups of $\Z_6$: $(e)$, $\Z_2$, $\Z_3$, and $\Z_6$. One has:
\begin{eqnarray*}
(1-t)^{-[\Z_6/(e)]}&=&\frac{1}{1-t^6}[1]+\frac{t^3}{(1-t^3)(1-t^6)}[\Z_6/\Z_3]\\
&+&\frac{t^2}{(1-t^2)^2(1-t^6)}[\Z_6/\Z_2]\\
&+&\frac{t(1+t+2t^2+2t^3+3t^4+2t^6+t^7)}{(1-t^2)^2(1-t^3)(1-t^6)(1-t)^2}[\Z_6/(e)],\\
(1-t)^{-[\Z_6/\Z_2]}&=&\frac{1}{1-t^3}[1]+\frac{t}{(1-t^3)(1-t)^2}[\Z_6/\Z_2],\\
(1-t)^{-[\Z_6/\Z_3]}&=&\frac{1}{1-t^2}[1]+\frac{t}{(1-t)(1-t^2)}[\Z_6/\Z_3].
\end{eqnarray*}

\item Let $G$ be the group $S_3$ of permutations on three elements. There are four conjugacy classes of subgroups of $G$: $(e)$, $\Z_2$, $\Z_3$, and $S_3$ (the class $\Z_2$ consists of three subgroups). One has: 

\begin{eqnarray*}
(1-t)^{-[S_3/(e)]}&=&\frac{1}{1-t^6}[1]+\frac{t^3}{(1-t^3)(1-t^6)}[S_3/\Z_3]\\
&+&\frac{3t^2}{(1-t^2)^2(1-t^6)}[S_3/\Z_2]\\
&+&\frac{t(1+4t^2+t^3+4t^4-2t^5+3t^6+t^7)}{(1-t^2)^2(1-t^3)(1-t^6)(1-t)^2}[S_3/(e)],\\
(1-t)^{-[S_3/\Z_2]}&=&\frac{1}{1-t^3}[1]+\frac{t(1+2t)}{(1-t^3)(1-t)^2}[S_3/\Z_2]\\
&+&\frac{t^3}{(1-t^2)(1-t)(1-t^2)}[S_3/(e)],\\
(1-t)^{-[S_3/\Z_3]}&=&\frac{1}{1-t^2}[1]+\frac{t}{(1-t)(1-t^2)}[S_3/\Z_3].
\end{eqnarray*}

\end{enumerate}
\end{examples}

\bigskip

There is a natural homomorphism from the Grothendieck ring $K_0(\fgs)$ to the ring $R(G)$ of virtual representations of the group $G$ which sends the class $[G/H]$ of a $G$-set  to the representation $i^G_H[1_H]$ induced from the trivial one-dimensional representation $1_H$ of the subgroup $H$.
This homomorphism is a homomorphism of $\lambda$-rings (\cite{knutsen}).

\section{An equivariant Euler characteristic of a $G$-variety}

Let the group $G$ act on a variety $X$. There exists a representation of $X$ as a finite cell complex such that the action of an element $g\in G$ respects the cell decomposition and, moreover, if $g~\sigma=\sigma$ for a cell $\sigma$, then $g_{|_\sigma}=id$. The set $C_n$ of $n$-cells is a finite $G$-set.

\begin{definition}
 The {equivariant Euler characteristic} $\chi_G(X)$ is the alternative sum $\sum_{n\geq 0}
(-1)^n [C_n]\in K_0(\fgs)$ of the classes of the $G$-sets $C_n$.
\end{definition}

\medskip

One has the following formula for the equivariant Euler characteristic $\chi_G(X)$
(which in particular shows that it does not depend on the cell decomposition). For a conjugacy class $h\in {\rm consub}(G)$ of subgroups of the group $G$, let $X_h:=\{x\in X:\, G_x\in h\}$, where $G_x:=\{g\in G:\, g\cdot x=x \}$ is the isotropy group of the point $x$.
Then one has 
\begin{equation}\label{eqeu}
\chi_G(X)=\sum_{h\in {\rm consub}(G)}\frac{\chi(X_h)\,|H|}{|G|}[G/H]=\sum_{h\in {\rm consub}(G)}\chi(X_h/G)[G/H],
\end{equation}
where $H$ is a representative of the conjugacy class $h$.

\begin{remark}
Here the Euler characteristic $\chi(\cdot)$ is the additive one, i.e. the alternative sum of ranks of the cohomology groups with compact support.
\end{remark}

The natural homomorphism from the Grothendieck ring $K_0(\fgs)$ to the ring $R(G)$ of virtual representations of the group $G$ sends the equivariant Euler characteristic $\chi_G(X)$ to the one defined by Wall in \cite{CTC}. Moreover, since this homomorphism is, generally speaking, not injective, the equivariant Euler characteristic as an element in $K_0(\fgs)$ is a somewhat finer invariant than the one as an element of the ring $R(G)$.

\section{Equivariant Lefschetz number and an e\-qui\-va\-riant zeta function of a map}

A concept of the {equivariant Lefschetz number} or {class} was introduced in \cite{LR} (see also \cite{GW}).

The {equivariant Lefschetz number} can be defined for a $G$-invariant map $\varphi$
from a finite $G$-$CW$-complex $X$ to itself and also for an equivariant map $\varphi$ from a pair $(X,Y)$ of finite $G$-$CW$-complexes to itself. In the last case the definition takes into account a $G$-equivariant chain complex of the pair $(X,Y)$ and the action of the map $\varphi_*$ on it. The latter means that the equivariant Lefschetz number can be defined also for a $G$-equivariant proper map from the space $X\setminus Y$ (which is not a finite $CW$-complex) to itself (by considering the one-point compactification of $X\setminus Y$).

The equivariant Lefschetz number $\Lambda^G(\varphi)$ of a $G$-map $\varphi:X\to X$ is an element of the Grothendieck ring $K_0(\fgs)$ of finite $G$-sets. If $X$ is a $G$-manifold, the equivariant Lefschetz number counts fixed points of a generic $G$-equivariant perturbation of the map $\varphi$ (which is a union of $G$-orbits) with appropiate signs.

We shall essentially use the following properties of the equivariant Lefschetz number:

\begin{enumerate}
\item Under the natural homomorphism from the Grothendieck ring $K_0(\fgs)$ to the ring $\Z$ of integers 
(the class of a $G$-set is sent to the number of points in it)
the equivariant Lefschetz number of a $G$-equivariant map $\varphi$ reduces to the usual Lefschetz number
$\Lambda(\varphi)$  of the map $\varphi$.

\item {\bf Additivity}. If $Y$ is a closed $G$-subspace of $X$ and a $G$-equivariant map  $\varphi:X\to X$ preserves $Y$ and $X\setminus Y$ (i.e. it sends each of these spaces to itself)
then 
$$
\Lambda^G(\varphi)=\Lambda^G(\varphi_{|_Y})+\Lambda^G(\varphi_{|_{X\setminus Y}}).
$$

\item If a $G$-equivariant map  $\varphi:X\to X$ has no fixed points then $\Lambda^G(\varphi)=0$.

\item $\Lambda^G(id_X)=\chi_G(X)$.

\item If $X=X_h=\{x\in X:\,G_x\in h\}$ for a conjugacy class $h$ of subgroups of the group $G$, then 
$$
\Lambda^G(\varphi)=\frac{\Lambda(\varphi)\,|H|}{|G|}[G/H],
$$
where $\Lambda(\varphi)$ is the usual Lefschetz number of the map $\varphi$.

\end{enumerate}

\medskip

In the usual (non-equivariant, i.e. $G=(e)$) situation, the zeta function of a map $\varphi:X\to X$ is the rational function in the variable $t$ which, as a series in $t$, is given by the formula
\begin{equation*}
\zeta_\varphi(t)=\prod_{m\geq 1} (1-t^m)^{-\frac{s_m}{m}},
\end{equation*}
where $s_m$ are the numbers
(in fact integers) defined by the formulae:

\begin{equation*}
\Lambda(\varphi^m)=\sum_{i|m}s_i.
\end{equation*}

The degree $\deg \zeta_\varphi$ of the zeta function of the map $\varphi$ (i.e. the degree of the numerator minus the degree of the denominator) is equal to the Euler characteristic of the variety $X$. This inspires the following definition:

\begin{definition}
The {equivariant zeta function} of a $G$-equivariant map $\varphi:X\to X$ is the series
$\zeta_\varphi^G(t)\in 1+t(K_0(\fgs)\otimes \Q)[[t]]$ defined by
\begin{equation}\label{eq2}
\zeta_\varphi^G(t)=\prod_{m\geq 1} (1-t^m)^{-\frac{s_m^{{{G}}}}{m}}, 
\end{equation}
where the virtual finite $G$-sets $s_m^G\in K_0(\fgs)$ are defined by the following recurrence formula:
\begin{equation}\label{eq3}
\Lambda^G(\varphi^m)=\sum_{i|m}s_i^G.
\end{equation}
and $(1-t^m)^{-\frac{s_m^G}{m}}$ is the series from $1+t(K_0(\fgs)\otimes \Q)[[t]]$
defined in Section~\ref{uno}.
\end{definition}

\begin{remark}
\begin{enumerate}
 \item Generally speaking, coefficients of the series $\zeta_\varphi^G(t)$
 do not belong to the Grothendieck ring $K_0(\fgs)$,
i.e. this series  may have ``rational'' coefficients.
\item Applying the natural homomorphism from the Grothendieck ring $K_0(\fgs)$
to the ring $R(G)$ of representations of the group $G$ (see Section \ref{uno}) one gets a version of the zeta function as an element of $1+t(R(G)\otimes \Q)[[t]]$.
\end{enumerate}

\end{remark}

\section{The A'Campo type formula for the equivariant zeta function of the monodromy}

Let $(V,0)$ be a germ of a purely $n$-dimensional complex analytic variety with an action of the group $G$ and let $f:(V,0)\to (\C,0)$ be the germ of a $G$-invariant analytic function such that
$\mbox{\rm Sing} \,V\subset f^{-1}\{0\}$. Let $M_f$ be the Milnor fibre of the germ $f$ at the origin: $M_f=\{x\in V: f(x)=\varepsilon, \Vert x \Vert\leq \delta\,\}$ with $0<\vert \varepsilon \vert \ll \delta $ small enough (for this definition we assume $(V,0)$  to be embedded in a certain affine space $(\C^N,0)$). The group $G$ acts on the Milnor fibre $M_f$ because the function germ $f$ is $G$-invariant.

Let $\pi:(X,\DD)\to (V,0)$ be a $G$-equivariant resolution of the germ $f$, i.e. a proper
$G$-equivariant map from a $G$-manifold $X$ to $V$ such that $\pi$ is an isomorphism outside the zero level set  $f^{-1}\{0\}$ of the function $f$ and, in a neighbourhood of any point $p$ of the total transform 
$E_0:=\pi^{-1}(f^{-1}\{0\})$ of the zero-level set $f^{-1}\{0\}$ of the function $f$, there exists a local system of coordinates $z_1,\ldots,z_n$
(centred at the point $p$) in which one has $f\circ \pi(z_1,\ldots,z_n)=z_1^{m_1}z_2^{m_2}\cdots z_n^{m_n}$,
with non-negative integers $m_i$, $i=1,\ldots,n$ (this implies that the total transform $E_0$ of the zero level set $f^{-1}\{0\}$ is a normal crossing divisor on $X$).
Moreover, we assume that, for each point 
$p\in E_0$, the irreducible components of $E_0$ are invariant with respect to the isotropy group $G_p$ of the point $p.$

Let $S_m$, $m\geq 1$, be the set of points $p$ of the exceptional divisor $\DD:=\pi^{-1}\{0\}$ such that, 
in a neighbourhood of the point $p$, one has $f\circ \pi(z_1,\ldots,z_n)=z_1^{m}$.
For $p\in S_m$, let $G_p$ be the isotropy group of the point $p$: $G_p=\{g\in G:\,gp=p\}$.
The group $G_p$ acts on the smooth germ $(X,p)$ preserving the exceptional divisor $\DD$ locally
given by $z_1=0$. This implies that, in a neighbourhood of the point $p$, one can suppose $G_p$ to act by linear transformations in coordinates $z_1,\ldots,z_n$ preserving ``the normal slice''
$z_2=\ldots=z_n=0$. This way one gets a linear representation of the group $G_p$ on this 
normal slice. 
Let ${\widehat G}_p$ be the kernel of this representation. One can see that $G_p/{\widehat G}_p$
is a cyclic group the order of which divides $m$.

 Let $S_{m,H,{\widehat H}}$, $({\widehat H}\subseteq H)$, be the set of points $p\in S_m$
such that the pair $(H,{\widehat H})$ is conjugate to the pair $ (G_p,{\widehat G}_p)$, i.e.
for the same element $g\in G$ one has  $G_p=gHg^{-1}$ and ${\widehat G}_p=g{\widehat H}g^{-1}$.

\begin{theorem*}\label{acampo}
\begin{eqnarray}
\zeta^G_f(t)&=&\prod_{m\geq 1,(h,{\widehat h})}(1-t^m)^{-\frac{|{\widehat H}|\chi(S_{m,H,{\widehat H}})}{|G|}[G/{\widehat H}]} \nonumber \\
&=&\prod_{m\geq 1,(h,{\widehat h})}(1-t^m)^{-\frac{|{\widehat H}|\chi(S_{m,H,{\widehat H}}/G)}{|H|}[G/{\widehat H}]}\,,\label{eq4}
\end{eqnarray}
where the product is over all conjugacy classes $(h,{\widehat h})$ of pairs of subgroups of the group $G$, the pair $(H,{\widehat H})$ is a representative of the conjugacy class $(h,{\widehat h})$.
\end{theorem*}

\begin{proof}
Just like in the non-equivariant case (see \cite{AC},\cite{clemens}) one can construct a $G$-equivariant retraction of a neighbourhood of the total transform $E_0$ of the zero level set $f^{-1}\{0\}$ to $E_0$ itself such that, outside of a  neighbourhood, in $E_0$, of all intersections of at least two irreducible components of $E_0$, i.e. where in local coordinates  $f\circ \pi(z_1,\ldots,z_n)=z_1^{m},\, m\geq 1$, the retraction sends a point $(z_1,z_2,\ldots,z_n)$ to $(0,z_2,\ldots,z_n)$.

Moreover, one can construct a monodromy transformation $\varphi$ which commutes with the retraction and, in a neighbourhood of a point where $f\circ \pi(z_1,\ldots,z_n)=z_1^{m}$, it sends
a point
$(z_1,z_2,\ldots,z_n)\in M_f$ to $(\exp(\frac{2\pi i}{m})z_1,z_2,\ldots,z_n)\in M_f$, $M_f$ being the Milnor fibre.

From the additivity property 2 above it follows that the equivariant Lefschetz number $\Lambda^G(\varphi^k)$ is equal to the sum of the equivariant Lefschetz numbers for the map $\varphi^k$ in a neighbourhood of intersections of the irreducible components of $E_0$ and also on the preimages of the strata $S_{m,H,{\widehat H}}$.

One can show that the equivariant Lefschetz number of the restriction of the power $\varphi^k$
of the monodromy transformation $\varphi$ to a neighbourhood of a $G$-orbit in the  intersection of irreducible components of $E_0$ is equal to zero. This follows from the fact that the Milnor fibre $M_f$ in a neighbourhood of such an orbit can be fibred by circles and the monodromy transformation can be supposed to preserve these fibration. Then the fact that the equivariant Lefschetz number is equal to zero follows from the fact that the Euler characteristic of the circle is equal to zero. Additivity of the equivariant Lefschetz number
implies that the equivariant Lefschetz number of the restriction of the power $\varphi^k$
of the monodromy transformation $\varphi$ to a neighbourhood of the union of intersections of irreducible components of $E_0$ is equal to zero.

The equivariant Lefschetz number of the restriction of the power $\varphi^k$
of the monodromy transformation to the preimage of a stratum $S_{m,H,{\widehat H}}$ can be computed using Property~5 of the equivariant Lefschetz number. This implies that:
\begin{equation*}
\Lambda^G(\varphi^k)=\sum_{m|k,H,{\widehat H}}\frac{m}{|H|/|{\widehat H}|}\chi(S_{m,H,{\widehat H}}/G)[G/{\widehat H}]=\sum_{m|k,H,{\widehat H}}\frac{m |{\widehat H}|}{|G|}\chi(S_{m,H,{\widehat H}})[G/{\widehat H}].
\end{equation*}
This implies that the classes $s_m^G$ defined by (\ref{eq3}) are equal to 
\begin{equation*}
s_m^G=\sum_{H,{\widehat H}}\frac{|{\widehat H}|}{|G|}\chi(S_{m,H,{\widehat H}})[G/{\widehat H}].
\end{equation*}
(in particular, there are only finitely many classes $s_m^G$ different from $0$). Now the definition $(\ref{eq2})$ implies the equation (\ref{eq4}).
\end{proof}

\begin{remark}
One can see that the degree of the zeta function $(\ref{eq4})$ understood as the (finite) sum
\begin{equation*}
\sum_{m,H,{\widehat H}}\frac{m |{\widehat H}|}{|G|}\chi(S_{m,H,{\widehat H}})[G/{\widehat H}]
\end{equation*}
is equal to  the equivariant Euler characteristic $\chi^G(M_f)$ of the Milnor fibre $M_f$.
\end{remark}

\begin{examples}
\begin{enumerate}
\item Let $f:(\C^2,0)\to (\C,0)$ be the germ of the analytic map defined by $f(x,y)=x^3-y^2$. Let the cyclic group $\Z_6=<\exp(2i\pi/6)>\subset \C^*$ act on $\C^2$ by $\lambda\cdot (x,y)=(\lambda^2x,\lambda^3y)$. The minimal embedded resolution of the curve $\{f=0\}$ (by three blowing-ups) is an equivariant one.
There are two strata $S_{2,(e),(e)}$ and $S_{3,(\Z_2),(e)}$ with the Euler characteristic equal
to zero. The strata $S_{2,(\Z_6),(\Z_3)}$ and $S_{3,(\Z_6),(\Z_2)}$ consist of one point each.
The stratum $S_{6,(\Z_6),(e)}$ consists of one component of the exceptional divisor (isomorphic
to the projective line $\C\P^1$) without three points. Its Euler characteristic is equal to $-1$. Now Theorem gives
 \begin{equation*}
\zeta^{\Z_6}_f(t)=(1-t^2)^{-\frac{1}{2}[\Z_6/\Z_3]}(1-t^3)^{-\frac{1}{3}[\Z_6/\Z_2]}(1-t^6)^{\frac{1}{6}[\Z_6/(e)]}.
\end{equation*}

\item Let $V=\{(x,y,z)\in \C^3:\, x+y+z=0\}$ and let $f$ be the function on $V$ 
Defined by $f(x,y,z)=x^2+y^2+z^2$.
Consider the natural action of the group $S_3$ of permutations on three elements on $V$ ($f$ is an $S_3$-invariant function). Blowing-up the origin, one gets a resolution of the function $f$. There are three strata of the
form $S_{m,H,{\widehat H}}$. The stratum $S_{2,(\Z_2),(\Z_2)}$ consists of three points $(1:1:-2)$,
$(1:-2:1)$, and $(-2:1:1)$ permuted by the group $S_3$ in the natural way. The stratum
$S_{2,(\Z_2),(e)}$ consists of three points $(1:-1:0)$, $(1:0:-1)$, and $(0:1:-1)$. The stratum
$S_{2,(e),(e)}$ coincides with the exceptional divisor $\C\P^1$ of the resolution without 8 points.
Its Euler characteristic is equal to $-6$. Now Theorem gives
 \begin{equation*}
\zeta^{S_3}_f(t)=(1-t^2)^{-[S_3/\Z_2]}(1-t^2)^{\frac{1}{2}[S_3/(e)]}.
\end{equation*}
\end{enumerate}
\end{examples}


\begin{thebibliography}{15}

\bibitem{AC} N.~A'Campo. {La fonction z\^eta d'une monodromie}. 
Comment. Math. Helv.  50  (1975), 233--248.

\bibitem{clemens} C.H.~Clemens. 
{Picard-Lefschetz theorem for families of nonsingular algebraic varieties acquiring ordinary singularities}.
Trans. Amer. Math. Soc. 136 (1969), 93--108. 

\bibitem{GW} D.L.~Goncalves,  J.~Weber.
{Axioms for the equivariant Lefschetz number and for the Reidemeister trace}.
Preprints of the Max-Planck-Institut f\"ur Mathematik, MPIM2006-123.

\bibitem{GLM1} S.M.~Gusein-Zade, I.~Luengo, A.~Melle-Hern\'andez.
{A power structure over the Grothendieck ring of varieties}.
Math. Res. Lett. 11 (2004), 49--57.

\bibitem{knutsen} D.~Knutson.  {$\lambda $-rings and the representation theory of the symmetric group}. Lecture Notes in Mathematics, Vol. 308. Springer-Verlag, Berlin-New York, 1973. iv+203 pp. 

\bibitem{LR} W.~L\"uck, J.~Rosenberg. {The equivariant Lefschetz fixed point theorem for proper cocompact $G$-manifolds}. In: High-dimensional manifold topology,  322--361, World Sci. Publ., River Edge, NJ, 2003. 

\bibitem{CTC} C.T.C.~Wall. {A note on symmetry of singularities}. Bull. London Math. Soc. 12 (1980), no.3, 169--175.
\end{thebibliography}
\end{document}